\documentclass[a4paper,11pt]{article}
\usepackage{hyperref}
\usepackage{amsmath,amssymb,amsfonts}  
\usepackage{graphics}
\title{Extensions of endomorphisms of $C(X)$} 
\author{J. F. Feinstein and T. J. Oliver} 
\date{}

\def\Z{\mathbb{Z}}
\def\T{\mathbb{T}}
\def\End{\textrm{End}}
\begin{document} 
\maketitle 
\begin{abstract} 
For a compact space $X$ we consider extending endomorphisms of the 
algebra $C(X)$ to be endomorphisms of Arens-Hoffman and Cole
extensions of $C(X)$. Given a non-linear,
monic polynomial $p\in C(X)[t]$, with $C(X)[t]/pC(X)[t]$ 
semi-simple, we show that if an endomorphism of $C(X)$ extends to the 
Arens-Hoffman extension with respect to $p$ then it also extends to the simple Cole 
extension with respect to $p$. We show that the converse to this is false. For 
locally connected, metric $X$ we characterize the algebraically closed $C(X)$, 
in terms of the extendability of endomorphisms to
Arens-Hoffman and to simple Cole extensions. 
\end{abstract} 

\section{Introduction} 
The problem of classifying the endomorphisms of commutative Banach algebras and their spectra
has been extensively studied (see, for example, \cite{FeiKam,GalGamLin,Kam,KamSch,Klein}). 
By contrast, the extendability problem for endomorphisms 
of commutative Banach algebras appears to be 
relatively untouched. In this note we concentrate on two different types of 
extension of the algebra $C(X)$ of continuous, complex-valued functions on 
$X$ 
obtained by adjoining a root of some polynomial $p \in C(X)[t]$: 
the Arens-Hoffman extension and the 
\lq simple' Cole extension. 
For a thorough treatment of these algebra extensions and how they relate to 
each 
other we refer the reader to \cite{Daw} (see also \cite{AreHof,Cole}). 
Throughout this paper, all Cole extensions discussed will be assumed to be simple in the
above sense.

The question of whether or not an endomorphism of $C(X)$ extends to a given 
Cole extension may be viewed as purely topological. The question of 
whether or not an endomorphism of $C(X)$ extends to an Arens-Hoffman 
extension is more algebraic in nature. By asking these two questions 
simultaneously we gain insight into the relationships between the
associated topological and algebraic conditions.
We shall see, in the case when $C(X)[t]/pC(X)[t]$ is semi-simple, that 
if an endomorphism of $C(X)$ extends to the Arens-Hoffman extension with 
respect to $p$ then it also extends to the Cole extension with 
respect to $p$. We shall also give an example to show that
the converse is false.These results are, 
perhaps, 
surprising because the Arens-Hoffman extension may be viewed as a sub-algebra 
of the Cole extension. Finally, for locally connected, compact metric spaces $X$, 
we shall characterize when $C(X)$ is algebraically closed in terms of
the extendability of endomorphisms  to Cole or Arens-Hoffman extensions.

\section{Preliminary notation, definitions and results} 
Let $W$ and $X$ be topological spaces and $A$ be a commutative algebra. We 
denote the ring of polynomials with coefficients in $A$ by $A[t]$. The set of 
continuous functions from $W$ to $X$ is denoted by $C(W,X)$. The set of 
complex-valued, continuous functions on $X$ is denoted by $C(X)$ which, when 
$X$ is compact, is a commutative, unital, Banach algebra under pointwise 
addition and multiplication with the uniform norm $\|g\|:= \sup\{|g(x)|: x\in 
X\}$. For $g\in C(X)$ we denote the zero set of $g$ by $Z(g)$. 
For $q:= q_{0}+\ldots q_{m}t^{m}\in C(X)[t]$ and a homomorphism 
$\alpha:C(X)\rightarrow A$ we define $q^{(\alpha)}= 
\alpha(q_{0})+\ldots\alpha(q_{m})t^{m}\in A[t]$.
For $\beta\in C(W,X)$ we 
define $\beta^{*}: C(X)\rightarrow C(W)$ by $\beta^{*}(g)=g\circ\beta$ 
($g\in C(X)$).
 
In our terminology a {\em compact space} is a compact, Hausdorff topological space.
For a compact space $X$ and $x\in X$, $\widehat{x}$ denotes the evaluation
character on $C(X)$ at $x$: $\widehat{x}(g)=g(x)$. 
Given a  monic polynomial $p\in C(X)[t]$  of degree at least $2$
we define
$X_{p}= \{(x,\lambda)\in X\times\mathbb{C}: p^{(\widehat{x})}(\lambda)= 0\}$. 
We define $\pi_{X} : X_{p}\rightarrow X$ and 
$\pi_{p} : X_{p}\rightarrow \mathbb{C}$ 
to be the restrictions to $X_{p}$ of the two coordinate projections.
The discriminant 
$D_{p}\in C(X)$ of $p$ is defined by $D_{p}(x)= \prod_{(i< 
j)}(\lambda_{i}-\lambda_{j})^{2}$ where $p^{(\widehat{x})}= 
(t-\lambda_{1})\times\ldots (t-\lambda_{n})$. 
If $\textrm{int}(Z(D_{p}))= 
\emptyset$ then $D_{p}$ is not a zero-divisor in $C(X)$ and
$p$ is of minimal degree in $C(X)[t]\setminus \{0\}$ such that 
$p^{(\pi_{X}^{*})}(\pi_{p})= 0$.
In this case, the Arens-Hoffman 
extension $C(X)[t]/pC(X)[t]$ of $C(X)$ with respect to $p$ is semi-simple and
is isomorphic to the algebra
$C(X)_{p}:= \{q^{(\pi_{X}^{*})}(\pi_{p}) : q\in C(X)[t]\} \subseteq C(X_{p})$
 (\cite{AreHof}). 
Note that $C(X)_{p}$ is the sub-algebra of $C(X_{p})$ 
generated by $\{\pi_{p}\}\cup\pi_{X}^{*}(C(X))$.
In this setting, the Cole extension of $C(X)$ with respect to $p$ is simply
$C(X_{p})$ and this is equal to the uniform closure of  $C(X)_{p}$ (\cite{Cole,Daw}).
For convenience, we shall describe monic polynomials $p$ of degree at least two and with 
$\textrm{int}(Z(D_{p}))=\emptyset$ as {\it admissible} polynomials.

Throughout this note, all algebras considered will be commutative, unital, complex algebras.
If $A$ is an algebra we denote the set of unital endomorphisms of $A$ by 
$\End(A)$. Let $B$ be an algebra and $\gamma: A\rightarrow B$ be an 
injective homomorphism. Clearly $\gamma(A)$ is a sub-algebra of $B$, and on $\gamma(A)$
we may define an inverse to $\gamma$ which we denote by $\gamma^{-1}$.
For $S\in \End(A)$ we see that $\gamma S\gamma^{-1}\in 
\End(\gamma(A))$. For $T\in \End(B)$ we see that 
$T|_{\gamma(A)}= \gamma S\gamma^{-1}$ if and only if $T\gamma= 
\gamma S$. In this case we say that \textit{$S$ extends to $B$ via $\gamma$}, 
and \textit{$T$ is an extension of $S$ via $\gamma$}. 

For a compact
space $X$, the map $\phi \mapsto \phi^*$ (with notation as above) is a bijection from
$C(X,X)$ to $\End(C(X))$. The inverse  map is given by 
$T \mapsto \phi_T$ where (regarding $X$ as a subset of $C(X)^*$ in the usual way)
$\phi_T=T^*|_X$. In this setting we say that $T$ is the endomorphism of  $C(X)$
\textit{induced by} $\phi_T$ and that $\phi_T$
is the self-map of $X$ {\it associated with} $T$.

We now give two observations which the reader may find useful to check before 
proceeding.
\begin{enumerate}
\item[(i)]
Let $X$ be a compact space, $p,q\in C(X)[t]$ be monic, and $f\in C(X_{p})$. 
Then $q^{(\pi_{X}^{*})}(f)= 0$ if and only if $(x,f(x,\lambda))\in X_{q}$ for 
all $(x,\lambda)\in X_{p}$. 
\item[(ii)]
Let $X$ be a compact space, $p\in C(X)[t]$ be monic, and $T\in 
\End(C(X))$. Then $(x,\lambda)\in X_{p^{(T)}}$ if and only if 
$(\phi_{T}(x),\lambda)\in X_{p}$. 
\end{enumerate}

\section{Extending endomorphisms} 
We begin our investigation with some lemmas concerning the extendability of
endomorphisms to the two types of algebra extension under consideration. 
 
Recall that we defined the term {\it admissible} for
polynomials in the previous section.\\

\textbf{Lemma 1}: Let $X$ be a compact space, let $p\in C(X)[t]$ be admissible, 
and let $T\in \End(C(X))$. 
\begin{enumerate}
\item[(i)]
There is an extension $T_{1}\in \End(C(X)_{p})$ of $T$ via 
$\pi_{X}^{*}$ if and only if there exists $f\in C(X)_{p}$ with 
$(x,f(x,\lambda))\in X_{p^{(T)}}$ for all $(x,\lambda)\in X_{p}$. 
\item[(ii)]
Given an $f$ satisfying the condition in (i),
we may find an extension $T_1$ such that 
$T_{1}(\pi_{p})=f$.\\
\end{enumerate}

\textit{Proof}: (i) Suppose that such an extension $T_{1}$ exists. 
Then $0= T_{1}(0)= 
T_{1}(p^{(\pi_{X}^{*})}(\pi_{p}))= p^{(\pi_{X}^{*}T)}(T_{1}(\pi_{p}))$. That 
is $p^{(\pi_{X}^{*}T)}$ has a root $T_{1}(\pi_{p})\in C(X)_{p}$. Set $f:= 
T_{1}(\pi_{p})$. Clearly $p^{(\pi_{X}^{*}T)}(f)= 0$, which means that 
$(x,f(x,\lambda))\in X_{p^{(T)}}$ for all $(x,\lambda)\in X_{p}$. 

Conversely, suppose 
that such an $f$ exists. Setting $T_{1}(\pi_{p}):= f$ and extending as an 
algebra homomorphism gives us a well defined endomorphism $T_{1}$. To check 
that $T_{1}$ is well defined is easy when we note that $p$ is monic and of 
minimal degree in $C(X)[t]\setminus \{0\}$ such that 
$p^{(\pi_{X}^{*})}(\pi_{p})= 0$. 

Statement (ii) follows from the last part of the argument above.
$\square$\\ 

\textbf{Lemma 2}: Let $X$ be a compact space, let $p\in C(X)[t]$ be 
admissible, and let $T\in \End(C(X))$. 
Then the following three statements are equivalent.
\begin{enumerate}
\item[(i)]
There exists an extension $T_{2}\in \End(C(X_{p}))$ of $T$ via 
$\pi_{X}^{*}$.
\item[(ii)]
There exists $\widetilde \phi \in C(X_p,X_p)$ such that $\pi_X  \circ \widetilde \phi = \phi_T \circ \pi_X$.
\item[(iii)]
There exists $f\in C(X_{p})$ with 
$(x,f(x,\lambda))\in X_{p^{(T)}}$ for all $(x,\lambda)\in X_{p}$. 
\end{enumerate}
Moreover, given such an $f$ satisfying the condition in (iii),
we may find an extension $T_2$ such that 
$T_{2}(\pi_{p})=f$.\\ 

\textit{Proof}: The equivalence of (i) and (ii) follows immediately from the
correspondence between self-maps and endomorphisms discussed above.

To see that (i) implies (iii),
suppose that such a $T_{2}$ exists. 
Then setting $f:= 
T_{2}(\pi_{p})\in C(X_{p})$ gives $p^{(\pi_{X}^{*}T)}(f)= 
p^{(T_{2}\pi_{X}^{*})}(T_{2}(\pi_{p}))= T_{2}(p^{(\pi_{X}^{*})}(\pi_{p}))= 0$, 
and so $(x,f(x,\lambda))\in X_{p^{(T)}}$ for all $(x,\lambda)\in X_{p}$. 

Finally we prove that (iii) implies (ii), and simultaneously prove the last part of the statement.
Suppose that such an $f$ exists as in (iii). Define $\widetilde \phi\in 
C(X_{p},X_{p})$ by $\widetilde \phi(x,\lambda):= (\phi_{T}(x),f(x,\lambda))$ for 
all $(x,\lambda)\in X_{p}$ . Set $T_2={(\widetilde \phi)}^*$.
Clearly (i) and (ii) hold and $T_{2}(\pi_{p})=f$. 
$\square$\\ 

Comparing Lemmas 1 and 2 we immediately obtain the following corollary.\\

\textbf{Corollary 3}: Let $X$ be a compact space, let $p\in C(X)[t]$ be admissible, and let
$T\in \End(C(X))$. If $T$ may be extended to
$C(X)_{p}$ via 
$\pi_{X}^{*}$ then $T$ may also be extended to  $C(X_{p})$ via $\pi_{X}^{*}$.\\

Another way to see this is to consider the operator norm of an extension of $T$ to
$C(X)_{p}$  when $C(X)_{p}$ is given the uniform norm rather than the Arens-Hoffman norm.
 
This result leads us to the following three questions.\\ 

\textbf{Question 1}: Do there exist a compact space $X$, admissible  $p\in C(X)[t]$, and 
$T\in \End(C(X))$ such that $T$ extends to $C(X)_{p}$ 
via $\pi_{X}^{*}$ but there exists an extension $T_{2}\in 
\End(C(X_{p}))$ of $T$ via $\pi_{X}^{*}$ with $T_{2}|_{C(X)_{p}}\notin 
\End(C(X)_{p})$?\\ 

\textbf{Question 2}: Do there exist a compact space $X$, admissible $p\in C(X)[t]$, and $T\in \End(C(X))$ 
such that $T$ extends to $C(X_{p})$ 
via $\pi_{X}^{*}$ but $T$ does not extend to 
$C(X)_{p}$ via $\pi_{X}^{*}$?\\ 

\textbf{Question 3}: Do there exist a compact space $X$, admissible $p\in C(X)[t]$, 
and $T\in \End(C(X))$ such that $T$ does not extend to 
$C(X_{p})$ via $\pi_{X}^{*}$?\\ 

In order to answer Question 1 we need the following lemma.\\ 

\textbf{Lemma 4}: Let $X$ be a compact space, let $p\in C(X)[t]$ be 
admissible, $f\in C(X)_{p}$, and $Y\subseteq X$ be a continuum. If $y_{0}$ is an 
interior point of $Y$, $\lambda_{1},\lambda_{2}\in C(Y)$ are such that 
$\lambda_{1}(y)\neq \lambda_{2}(y)$ for all $y\in Y\setminus \{y_{0}\}$, and 
$(y,\lambda_{1}(y)),(y,\lambda_{2}(y))\in X_{p}$ for all $y\in Y$ then 
$\lim_{y\rightarrow 
y_{0}}\left(\frac{f(y,\lambda_{1}(y))-f(y,\lambda_{2}(y))}{\lambda_{1}(y)-\lambda_{2}(y)}\right)$ 
exists and is finite.\\ 

\textit{Proof}: Choose $q:= q_{0}+\ldots q_{m}t^{m}\in C(X)[t]$ such that 
$q^{(\pi_{X}^{*})}(\pi_{p})= f$. For all $y\in Y\setminus\{y_{0}\}$ we have 
\[
\frac{f(y,\lambda_{1}(y))-f(y,\lambda_{2}(y))}{\lambda_{1}(y)-\lambda_{2}(y)}
= \sum_{k=1}^m h_k(y),
\]
where
$h_k(y)=q_{k}(y)(\lambda_{1}(y)^{k-1}+\lambda_{1}(y)^{k-2}\lambda_{2}(y)+\ldots 
\lambda_{2}(y)^{k-1}).
$
It follows immediately that the required limit exists and is finite. 
$\square$\\ 

We now answer Question 1 in the affirmative by way of an example.\\ 

\textbf{Example 1}: Let $X:= [0,1]$ and define $r\in C(X)$ by $r(x):= 
(3x-1)(3x-2)^{2}$ for all $x\in X$. The properties of $r$ that will be useful 
to us are that $r(\frac{1}{3})= r(\frac{2}{3})= 0$, $r'(\frac{1}{3})\neq 0$, 
and $r'(\frac{2}{3})= 0$. Now set $p(t):= (t-r)(t+r)\in C(X)[t]$ and let 
$T$ be the automorphism of
$C(X)$ induced by the self-homeomorphism $\phi(x)=1-x$ of $X$.

Setting $f\in C(X_{p})$ to be $f(x,\lambda):= r(1-x)$ for $\lambda= r(x)$, and 
$f(x,\lambda):= -r(1-x)$ for $\lambda= -r(x)$ we have that 
$(x,f(x,\lambda))\in X_{p^{(T)}}$ for all $(x,\lambda)\in X_{p}$. So, by Lemma 
2, there exists an extension $T_{2}\in \End(C(X_{p}))$ of $T$ via 
$\pi_{X}^{*}$ with $T_{2}(\pi_{p})= f$. We see that $f\notin C(X)_{p}$ 
otherwise we would have, by Lemma 4, that $\lim_{x\rightarrow 
\frac{2}{3}}\left(\frac{f(x,r(x))-f(x,-r(x))}{r(x)-(-r(x))}\right)= 
\lim_{x\rightarrow \frac{2}{3}}\left(\frac{r(1-x)}{r(x)}\right)$ exists and is finite.
However, 
the aforementioned properties of $r$ show that this is not the case. Thus, by 
Lemma 1, this particular $T_{2}$ does not restrict to give an endomorphism of 
$C(X)_{p}$.

Setting $f\in C(X_{p})$ to be $f(x,\lambda):= r(1-x)$ for all $(x,\lambda)\in 
X_{p}$ we still have that $(x,f(x,\lambda))\in X_{p^{(T)}}$ for all 
$(x,\lambda)\in X_{p}$, but we also have that $f\in C(X)_{p}$, since $f= 
\pi_{X}^{*}(T(r))$. So, by Lemma 1, $T$ extends to $C(X)_{p}$ via 
$\pi_{X}^{*}$. $\square$\\ 

In order to answer Question 2 we need another lemma.\\ 

\textbf{Definition}: Let $S\subset \T\times\mathbb{C}$, where 
$\T:= \{\lambda\in \mathbb{C}: |\lambda|= 1\}$. We call $S$ a 
\textit{$k$-strip with respect to $\mu$} if $k\in \mathbb{N}$ and $\mu: 
\T\rightarrow S$ is a homeomorphism such that $\mu(e^{i\theta})\in 
\{e^{ik\theta}\}\times\mathbb{C}$ for all $e^{i\theta}\in \T$.\\ 

\textbf{Lemma 5}: Let $a,b\in \mathbb{N}$, $A$ be an $a$-strip with respect to 
$\mu_{A}$, and $B$ be a $b$-strip with respect to $\mu_{B}$. There exists 
$g\in C(A,B)$ with $g((\{e^{i\theta}\}\times\mathbb{C})\cap A)\subset 
\{e^{i\theta}\}\times\mathbb{C}$ for all $e^{i\theta}\in \T$ if and 
only if $b$ divides $a$. In this case $g(A)= B$.\\ 

\textit{Proof}: Suppose that such a $g$ exists. Consider $\mu_{B}^{-1}\circ 
g\circ\mu_{A}\in C(\T,\T)$. Due to the properties of 
$\mu_{A},\mu_{B}$, and $g$ we see that $\mu_{B}^{-1}\circ 
g\circ\mu_{A}(e^{i\theta})= e^{i(\frac{a}{b}\theta+\theta_{0})}$ for all 
$\theta\in \mathbb{R}$, where $\mu_{B}^{-1}\circ g\circ\mu_{A}(e^{i0})= 
e^{i\theta_{0}}$. Now as $\mu_{B}^{-1}\circ g\circ\mu_{A}$ is continuous we 
must have that $\mu_{B}^{-1}\circ g\circ\mu_{A}(e^{i2\pi})= \mu_{B}^{-1}\circ 
g\circ\mu_{A}(e^{i\theta_{0}})$, which means that 
$e^{i(\frac{a}{b}2\pi+\theta_{0})}= e^{i\theta_{0}}$. Thus $b$ divides $a$. We 
now have that $\mu_{B}^{-1}\circ g\circ\mu_{A}(\T)= \T$ and 
so, as $\mu_{B}$ is bijective, $g(A)= B$. Suppose now that $b$ divides $a$. As 
$\mu_{A},\mu_{B}$ are bijective we may implicitly define $g$ by 
$\mu_{B}^{-1}\circ g\circ\mu_{A}(e^{i\theta}):= e^{i\frac{a}{b}\theta}$ for 
all $\theta\in \mathbb{R}$. It is clear that $g$ has the required properties. 
$\square$\\ 

We now give an example to show that Question 2  has an affirmative answer,
and that we may take $X$ to be $\T$.\\ 

\textbf{Example 2}: Define $\beta: [0,2\pi]\rightarrow \T$ by $\beta: 
\theta\mapsto e^{i\theta}$. We will define $p\in C(\T)[t]$ so that
$[0,2\pi]_{p^{(\beta^{*})}}$ has a certain desired structure. We shall choose 
$\lambda_{1},\ldots \lambda_{5}\in 
C([0,2\pi])$ and then choose an admissible $p\in C(\T)[t]$  so that
$[0,2\pi]_{p^{(\beta^{*})}}:= \bigcup_{\theta\in 
[0,2\pi]}(\{\theta\}\times \{\lambda_{1}(\theta),\ldots 
\lambda_{5}(\theta)\})$. 

We begin by choosing $\lambda_{1},\ldots \lambda_{5}\in 
C([0,2\pi])$ with $\lambda_{1}(2\pi)= \lambda_{2}(0)$, $\lambda_{2}(2\pi)= 
\lambda_{1}(0)$, $\lambda_{3}(2\pi)= \lambda_{4}(0)$, $\lambda_{4}(2\pi)= 
\lambda_{5}(0)$, $\lambda_{5}(2\pi)= \lambda_{3}(0)$, and such that 
$\lambda_{i}(\theta)= \lambda_{j}(\theta)$ for $i\neq j$ if and only if 
$\{i,j\}= \{2,5\}$ and $\theta= \pi$. When choosing $\lambda_{2}$ and 
$\lambda_{5}$ we insist in addition that, for $\theta \in [\pi-1,\pi+1] $, we have
$\lambda_{2}(\theta)= 
(\theta-\pi)^{2}$ and $\lambda_{5}(\theta)= -(\theta-\pi)^{2}$. 
It is now clear that we may choose an admissible quintic polynomial 
$p\in C(\T)[t]$ such that $p^{(\beta^{*})}(t)=\prod_{j=1}^5 (t-\lambda_j)$.
The resulting set $[0,2\pi]_{p^{(\beta^{*})}}$ is illustrated in 
Figure 1.

\begin{center}
\resizebox{11 cm}{8.75 cm}{\includegraphics{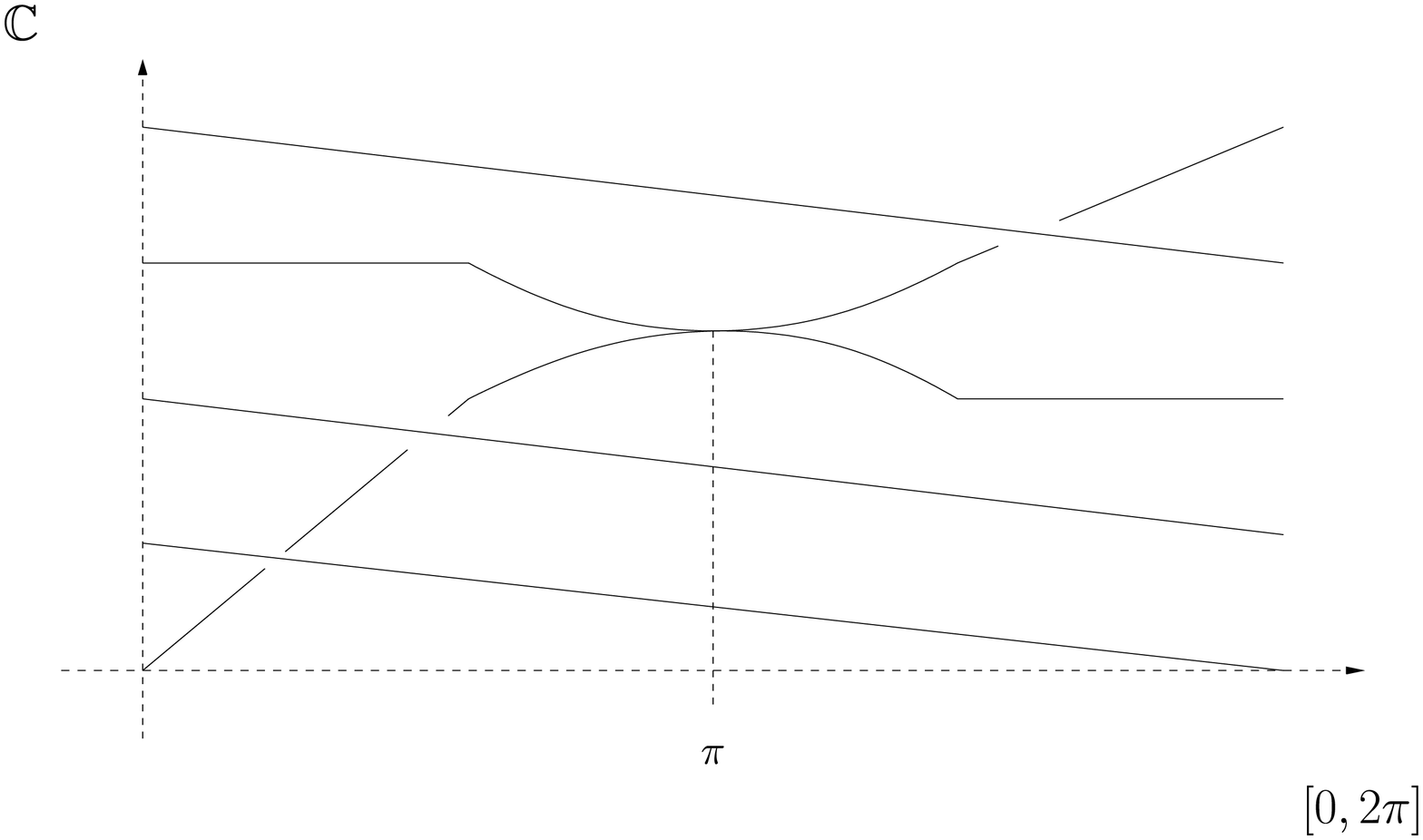}}\title{Figure 1}
\end{center}

Let $T$ be the endomorphism of $C(X))$ induced by the map $\phi$
where  $\phi(e^{i\theta}):= 
e^{i(\pi-\sqrt{\pi-\theta})}$ for $\theta\in [\pi-1,\pi]$, 
$\phi(e^{i\theta}):= e^{i(\pi+\sqrt{\theta-\pi})}$ for $\theta\in 
[\pi,\pi+1]$, and $\phi(e^{i\theta}):= e^{i\theta}$ otherwise. The set 
$[0,2\pi]_{p^{(\beta^{*}T)}}$ is illustrated in Figure 2. 

\begin{center}
\resizebox{11 cm}{8.75 cm}{\includegraphics{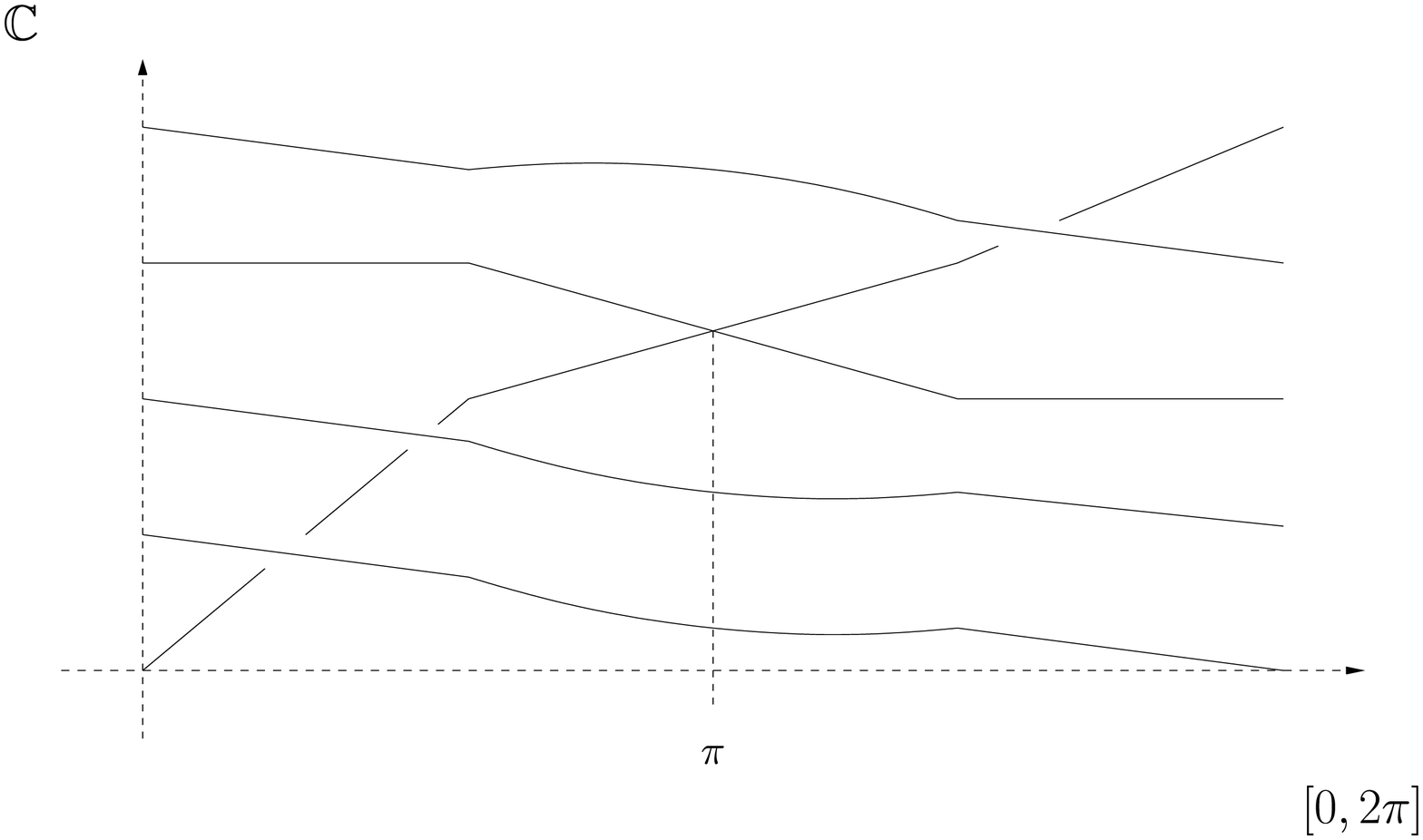}}\title{Figure 2}
\end{center}

The sets $\T_p$ and $\T_{p^{(T)}}$ 
are given by
\[
\T_p=\{(e^{i \theta},\lambda):(\theta,\lambda) \in [0,2\pi]_{p^{(\beta^*)}}\},
\]
\[
\T_{p^{(T)}}=\{(e^{i \theta},\lambda):(\theta,\lambda) \in [0,2\pi]_{p^{(\beta^* T)}}\}.
\]
Thus  each of the sets $\T_{p}$, $\T_{p^{(T)}}$ 
are the union of a 2-strip with a 3-strip. Denote the 2-strip in 
$\T_{p}$ by $R_{2}$, the 3-strip in $\T_{p}$ by $R_{3}$, the 
2-strip in $\T_{p^{(T)}}$ by $S_{2}$, and the 3-strip in 
$\T_{p^{(T)}}$ by $S_{3}$. In view of Lemma 5, setting $g(x,\lambda):= 
(x,f(x,\lambda))$ for all $(x,\lambda)\in \T_{p}$, we see that there 
is exactly one $f\in C(\T_{p})$ such that $(x,f(x,\lambda))\in 
\T_{p^{(T)}}$ for all $(x,\lambda)\in \T_{p}$. This $f$ is 
such that $g$ drops $R_{2}$ onto $S_{2}$, and $R_{3}$ onto $S_{3}$. Now by a 
similar argument to the one used in Example 1, considering 
$\lim_{\theta\rightarrow \pi}
\left(
\frac{f(e^{i\theta},\lambda_{2}(\theta))-f(e^{i\theta},\lambda_{5}(\theta))}{\lambda_{2}(\theta)-\lambda_{5}(\theta)}
\right)$, 
we see that $f\notin C(\T)_{p}$. So $T$ extends to $C(\T_{p})$ 
via $\pi_{\T}^{*}$ but $T$ does not extend to $C(\T)_{p}$ via 
$\pi_{\T}^{*}$. $\square$\\ 

Our third example shows that the answer Question 3 is also affirmative, and again 
we may take $X$ to be $\T$.\\ 

\textbf{Example 3}: Choose $p\in C(\T)[t]$ as in Example 2. Let $T\in 
\End(C(\T))$ be associated with $\phi_{T}\in 
C(\T,\T)$ where $\phi_{T}(e^{i\theta}):= e^{i(\theta+\pi)}$ 
for all $\theta\in \mathbb{R}$. Denote the 2-strip in $\T_{p}$ by 
$R_{2}$, and the 3-strip in $\T_{p}$ by $R_{3}$. Finally, denote the 
2-strip in $\T_{p^{(T)}}$ by $S_{2}$, and the 3-strip in 
$\T_{p^{(T)}}$ by $S_{3}$.

Suppose, for a contradiction, that there exists $f\in C(\T_{p})$ such 
that $(x,f(x,\lambda))\in \T_{p^{(T)}}$ for all $(x,\lambda)\in 
\T_{p}$. Define $g\in C(\T_{p},\T_{p^{(T)}})$ by 
$g(x,\lambda):= (x,f(x,\lambda))$ for all $(x,\lambda)\in \T_{p}$. We 
see, as in Example 2, that $g$ must drop $R_{2}$ onto $S_{2}$, and $R_{3}$ 
onto $S_{3}$. Hence $g$ is surjective. We note that $(\{-1\}\times 
\mathbb{C})\cap \T_{p}$ contains exactly four elements but 
$(\{-1\}\times \mathbb{C})\cap \T_{p^{(T)}}$ contains exactly five 
elements, as $\phi_{T}(-1)\neq -1$. This is a contradiction since $g$ is 
single valued. Hence no such $f\in C(\T_{p})$ may exist. Thus, by 
Lemma 2, $T$ does not extend to $C(\T_{p})$ via 
$\pi_{\T}^{*}$. $\square$\\

Careful consideration of the possible structures for $X_{p}$ and $X_{p^{(T)}}$
shows that, when $X=\T$, no admissible polynomial of degree less than $5$ 
can lead to a positive answer to Question 3. If, however, you move to higher-dimensional
$X$ the situation is different. The following simple example of this phenomenon 
was suggested to us by J.~W.~Barrett. Take $X=\T \times \T$, let $\phi$ be the self-homeomorphism of $X$
given by $\phi(z,w)=(w,z)$ and let $T$ be the automorphism of $C(X)$ induced by $\phi$.
Now let $F$ be the projection onto the first coordinate and
let $p$ be the quadratic polynomial $t^2-F$ (which is clearly admissible, and indeed has 
invertible discriminant).
Then it is easy to see that $T$ has no extension to $C(X_p)$ via $\pi_X^*$.

\section{Characterizing the algebraically closed $C(X)$}

An algebra $A$ is said to be {\it algebraically closed} if every non-constant, 
monic polynomial with coefficients in $A$ has a root in $A$.
In \cite{Cou}, Countryman characterized the first-countable compact spaces for which
$C(X)$ is algebraically closed in terms of certain hereditary topological conditions on 
$X$ which we describe below. He also showed that, for such $X$,
$C(X)$ is algebraically closed if and only if $C(X)=\{g^2:g\in C(X)\}$,
and gave an example to show that this latter equivalence fails if $X$ 
is not assumed to be first-countable.

For locally connected compact spaces, a different approach was taken by Hatori, Miura and Niijima
(\cite{HatMiu,MiuNii}).
They showed that, for such $X$, $C(X)$ is algebraically closed if and only if
the covering dimension of $X$ is at most $1$ and $H^1(X,\Z)$ is trivial,
where $H^1(X,\Z)$ is the first Cech cohomology group with
integer coefficients. 
(It is perhaps worth noting that the combination of the latter two conditions 
is also equivalent to the condition that
$\exp(C(X))$ be dense in $C(X)$.)
For locally connected compact spaces it is again 
sufficient to be able to factorize the quadratic monic polynomials:
$C(X)$ is algebraically closed if and only if $C(X)=\{g^2:g\in C(X)\}$. 

We now investigate the connection between $C(X)$ being algebraically closed 
and the extendability of endomorphisms to Arens-Hoffman 
extensions and Cole extensions. We begin with a lemma which immediately 
shows the connection in one direction.\\ 

\textbf{Lemma 6}: Let $X$ be a compact space, $p\in C(X)[t]$ be 
admissible, and $T\in \End(C(X))$. 
If $p^{(T)}(r)= 0$ for some $r\in C(X)$ then $T$ extends to $C(X)_{p}$ via 
$\pi_{X}^{*}$.\\ 

\textit{Proof}: Clearly $(x,r(x))\in X_{p^{(T)}}$ for all $(x,\lambda)\in 
X_{p}$. Set $f:= \pi_{X}^{*}(r)\in C(X)_{p}$ and apply Lemma 1. $\square$\\ 

We now introduce the topological conditions used by Countryman to characterize 
when $C(X)$ is algebraically closed.\\

\textbf{Definition}: Let $X$ be a compact space. We say that $X$ is 
\textit{hereditarily unicoherent} if every pair of continua in $X$ have 
connected or empty intersection. We say that $X$ is \textit{not almost locally 
connected due to $((C_{n}),(x_{n}),(y_{n}))$} if $(C_{n})$ is a sequence of 
pairwise disjoint continua in $X$ which are open in $\overline{\bigcup_{n\in 
\mathbb{N}}C_{n}}$ such that $x_{n},y_{n}\in C_{n}$ for all $n$, 
$x_{n}\rightarrow x_{0}\in X$, and $y_{n}\rightarrow y_{0}\in 
X\setminus\{x_{0}\}$. If no such $((C_{n}),(x_{n}),(y_{n}))$ exist then $X$ is
\textit{almost locally connected}.\\ 

We may now state Countryman's characterization \cite{Cou}: for a first-countable
compact space $X$, $C(X)$ is algebraically closed if and only if $X$ is
hereditarily unicoherent and almost locally connected. 

For the rest of this paper we take an arc to mean a homeomorphic copy of 
$[0,1]$.\\ 

\textbf{Lemma 7}: Every almost locally connected, metric continuum is arc 
connected.\\ 

\textit{Proof}: In \cite{Cou} it is shown that every sequentially compact, 
almost locally connected Hausdorff continuum is locally connected. Thus every 
almost locally connected, metric continuum is locally connected. Every locally 
connected, metric continuum is arc connected, by Theorem 3.15 of 
\cite{HocYou}. Therefore every almost locally connected, metric continuum is 
arc connected. $\square$\\ 

\textbf{Lemma 8}: Let $X$ be a locally connected, compact metric space. Then 
$C(X)$ is algebraically closed if and only if $X$ does not contain a 
homeomorphic copy of $\T$.\\ 

\textit{Proof}: We use Countryman's characterization, described above.
It is clear that if $X$ contains a homeomorphic copy of $\T$
then $X$ is not hereditarily unicoherent, and so $C(X)$ is not algebraically closed.
Conversely, suppose that $C(X)$ is not algebraically closed.

Consider first the case where $X$ is almost locally connected but not hereditarily 
unicoherent. We may choose continua $M,N\subset X$ such that $M\cap N= A\cup 
B$ where $A$ and $B$ are disjoint, non-empty compact subsets of $X$. 
Since $M$ and $N$ must also 
be almost locally connected we know, by Lemma 7, that $M$ and $N$ are arc 
connected. It is now easy to see that $M\cup N$ contains a homeomorphic copy 
of $\T$.

The remaining case is where $X$ is not almost locally connected due to some
$((C_{n}),(x_{n}),(y_{n}))$. Let $X_{0}$ be the component of $X$ containing 
$x_{0}$ and $y_{0}$. As $X$ is locally connected we see that $X_{0}$ is open. 
So there is an $n_{0}\in\mathbb{N}$ such that $C_{n}\subset X_{0}$ for $n\geq 
n_{0}$. Each $C_{n}$ is open in $\overline{\bigcup_{n\in \mathbb{N}}C_{n}}$, 
and $X_{0}$ is locally arc connected, by Lemma 3.29 of \cite{HocYou}. So for 
each $n\geq n_{0}$, and $x\in C_{n}$ there is a $d_{x}> 0$ such that 
$B_{d_{x}}(x)$ is arc connected, and $B_{2d_{x}}(x)\cap C_{m}= \emptyset$ for 
all $m\neq n$. As $C_{n}$ is compact we may cover $C_{n}$ with finitely many 
such $B_{d_{x}}(x)$. Hence we may assume that the $C_{n}$ are arcs from 
$x_{n}$ to $y_{n}$. We now have $x_{n}\rightarrow x_{0}$, $y_{n}\rightarrow 
y_{0}$, each $C_{n}$ is an arc from $x_{n}$ to $y_{n}$, and there exist 
disjoint, arc connected neighbourhoods of $x_{0}$ and $y_{0}$. It follows easily that 
$X_{0}$ 
contains a homeomorphic copy of $\T$. $\square$\\ 

\textbf{Theorem 9}: Let $X$ be a locally connected, compact metric space. Then 
the following are equivalent. 
\begin{enumerate}
\item[(i)]
$C(X)$ is algebraically closed. 
\item[(ii)]
For all closed $Y\subseteq X$, admissible $p\in C(Y)[t]$, 
and $T\in \End(C(Y))$ we have that 
$T$ extends to $C(Y)_{p}$ via $\pi_{Y}^{*}$. 
\item[(iii)]
For all closed $Y\subseteq X$, admissible $p\in C(Y)[t]$, and 
$T\in \End(C(Y))$ we have 
that $T$ extends to $C(Y_{p})$ via $\pi_{Y}^{*}$.\\ 
\end{enumerate}

\textit{Proof}: Suppose, first, that (i) does not hold. Then, by Lemma 8, 
there is a $Y\subseteq X$ and a homeomorphism $h:Y\rightarrow\T$. 
Choose $p\in C(\T)[t]$ and $T\in \End(C(\T))$ as in 
Example 3. Let $T_{h}\in \End(C(Y))$ be associated with 
$h^{-1}\circ\phi_{T}\circ h\in C(Y,Y)$. Since $T$ does not extend to 
$C(\T_{p})$ via $\pi_{\T}^{*}$ we have that $T_{h}$ does not 
extend to $C(Y_{p^{(h^{*})}})$ via $\pi_{Y}^{*}$. Hence (iii) does not hold, and so, by Corollary 3
(ii) does not hold. 

Suppose, now, that (i) holds. Then for all closed $Y\subseteq X$ we have that 
$C(Y)$ is algebraically closed. Thus for admissible $p\in C(Y)[t]$
and $T\in \End(C(Y))$ we have that 
$p^{(T)}$ has a root in $C(Y)$. Lemma 6 now shows that (ii) 
holds and hence, by Corollary 3, (iii) holds. 
$\square$\\ 

We conclude with a conjecture.\\ 

\textbf{Conjecture}: Let $X$ be a first countable, compact space. 
Then the following are equivalent.
\begin{enumerate}
\item[(i)]
$C(X)$ is algebraically closed. 
\item[(ii)]
For all admissible $p\in C(X)[t]$ and 
$T\in \End(C(X))$ we have that $T$ extends to $C(X)_{p}$ via 
$\pi_{X}^{*}$. 
\item[(iii)]
For all admissible $p\in C(X)[t]$
and $T\in \End(C(X))$ we have that $T$ extends to $C(X_{p})$ via 
$\pi_{X}^{*}$. 
\end{enumerate}

\newpage
{\sf  School of Mathematical Sciences

 University of Nottingham

 Nottingham NG7 2RD, England

 Joel.Feinstein@nottingham.ac.uk
 
 Thomas.Oliver@maths.nottingham.ac.uk

\vskip 0.3cm
2000 Mathematics Subject Classification: 46J10, 47B48
}
\vskip 0.3cm
\noindent
The second author would like to thank the EPSRC for providing support for
this research
 
\end{document}